\documentclass[journal]{IEEEtran}
\AtBeginDvi{} 

\usepackage[pdftex]{graphicx} 
\usepackage{framed}
\usepackage{bm}
\usepackage{comment}
\usepackage[verbose]{cite}
\usepackage{algorithm}
\usepackage{algorithmic}
\usepackage{slashbox}
\newtheorem{theorem}{Theorem}

\newtheorem{lemma}{Lemma}
\newtheorem{corollary}{Corollary}

\usepackage{latexsym}

\usepackage{ascmac}

\ifCLASSINFOpdf
\else
\fi
\usepackage{amsmath}
\begin{document}
%
\title{Optimal leader selection and demotion in leader-follower multi-agent systems}
%
%
%

\author{Kazuhiro~Sato
\thanks{K. Sato is with the School of Regional Innovation and Social Design Engineering,
 Kitami Institute of Technology,
 Hokkaido 090-8507, Japan,
email: ksato@mail.kitami-it.ac.jp}
}

\maketitle
\thispagestyle{empty}
\pagestyle{empty}

\begin{abstract}
We consider  leader-follower multi-agent systems that have many leaders, defined on any connected weighted undirected graphs, and address
the leader selection and demotion problems.
The leader selection problem is formulated as a minimization problem for the $H^2$ norm of the difference between the transfer functions of the original and new agent systems, under the assumption that the leader agents to be demoted are fixed. 
The leader demotion problem is that of finding optimal leader agents to be demoted, and is formulated
using the global optimal solution to the leader selection problem.
We prove that a global optimal solution to the leader selection problem is the set of the original leader agents except for those that are demoted to followers.
To this end, we relax the original problem into a differentiable problem.
Then, by calculating the gradient and Hessian of the objective function of the relaxed problem, we prove that the function is convex.
It is shown that zero points of the gradient are global optimal solutions to the leader selection problem, which  is a finite combinatorial optimization problem. 
Furthermore, we prove that any set of leader agents to be demoted subject to a fixed number of elements is a solution to the leader demotion problem.
By combining the solutions to the leader selection and demotion problems, we prove that
 if we choose new leader agents from the original ones except for those specified by the set of leader agents to be demoted,
then the relative $H^2$ error between the  transfer functions of the original and new agent systems is completely determined by the numbers of  original leader agents and  leader agents that are demoted to follower agents.
That is, we reveal that the relative $H^2$ error does not depend on the number of agents on the graph. 
Finally, we verify the solutions using a simple example.
\end{abstract}

\begin{IEEEkeywords}
Leader-follower multi-agent system, leader selection, leader demotion
\end{IEEEkeywords}

%
\IEEEpeerreviewmaketitle

\section{Introduction} \label{sec1}
%
%
%
%
\IEEEPARstart{L}{eader} selection is an important issue in leader-follower multi-agent systems \cite{mesbahi2010graph, oh2015survey, olfati2004consensus, rahmani2009controllability},
which include, for example, vehicle formation control \cite{ren2007distributed} and sensor networks \cite{yu2009distributed}.
This is because leader agents influence the dynamics of follower agents.
The study \cite{commault2013input} has investigated structural modifications of leader-follower multi-agent systems resulting from a leader selection, and mechanisms that lead to controllability.
The authors of \cite{clark2014supermodular} introduced an analytical approach for selecting leader agents to minimize the total mean-square error of the
follower agent states, using their desired value in a steady-state in the presence of noisy communication links.
The authors of \cite{clark2014minimizing}  studied a leader selection problem for minimizing convergence errors experienced by follower agents.
The study \cite{patterson2017optimal} has also addressed similar problems to \cite{clark2014supermodular, clark2014minimizing},
where the work in \cite{patterson2017optimal} was limited to one-dimensional networks, such as a path graph and ring graph, and
more efficient algorithms were provided for solving the problems.

\begin{figure}[t]
\begin{center}
\vspace{-45mm}
\includegraphics[scale=0.4]{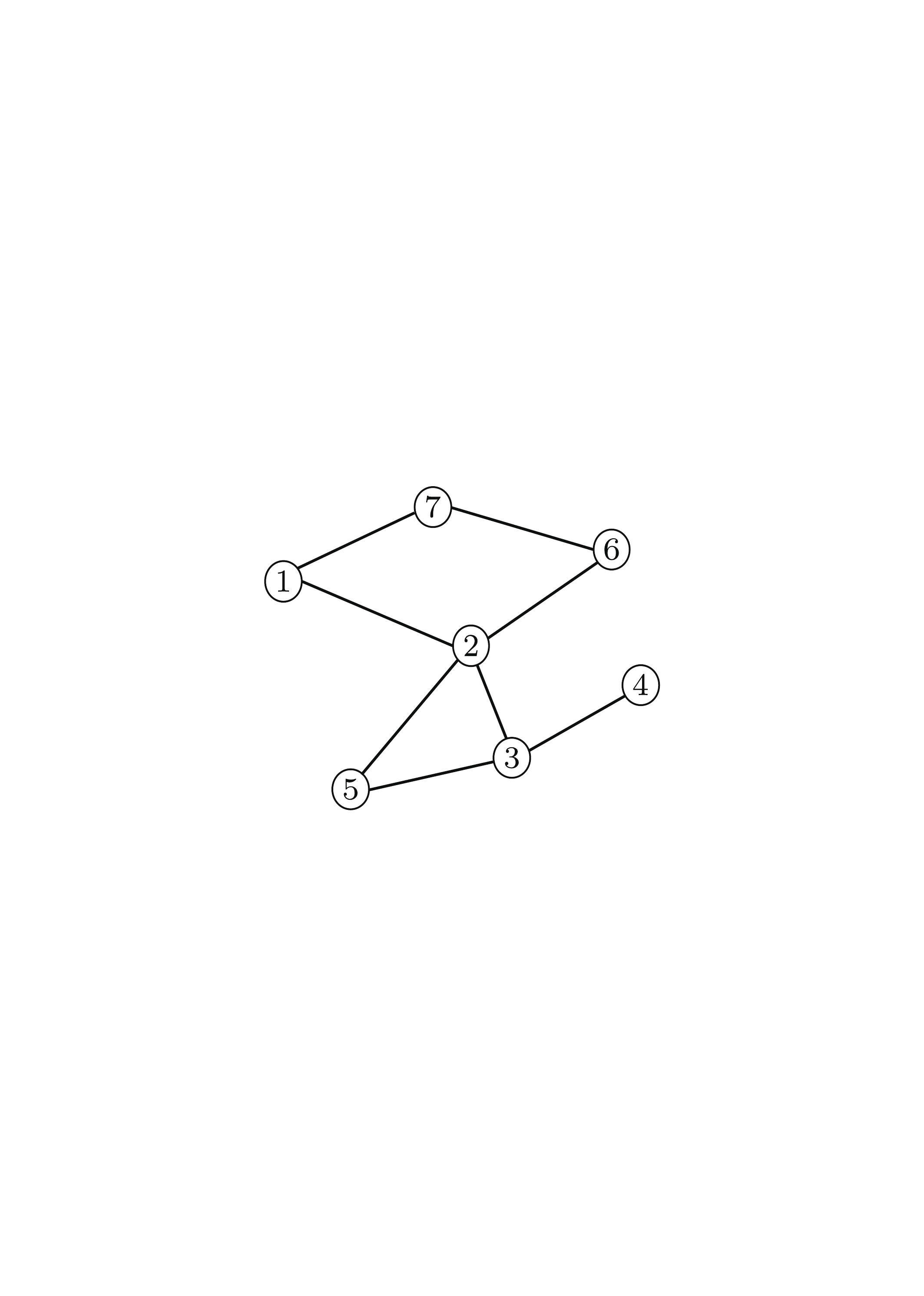}
\vspace{-45mm}
\end{center}
\caption{An example of connection relationship of multi-agent systems.} \label{graph_ex}
\end{figure}

In this paper, we consider leader-follower multi-agent systems that have many leaders, defined on any connected weighted undirected graphs.
Because these systems have many leaders, some leaders may not be important.
That is, the overall performance may not be affected, even if we demote some leaders to followers.
For example, suppose that $\{1,2,4\}$ in Fig.\,\ref{graph_ex} is a set of leader agents; i.e., $\{3, 5, 6, 7\}$ is a set of follower agents.
Furthermore, suppose that even if we demote leader agent $1$ to a follower,  the overall performance is almost unaffected.
Then, the multi-agent system with leaders $\{2,4\}$ and followers $\{1,3,5,6,7\}$ behaves like the original system with leaders $\{1,2,4\}$ and followers $\{3,5,6,7\}$.
However, there is a possibility that the dynamics of the multi-agent system with leaders $\{5, 6\}$ and followers $\{1, 2, 3, 4, 7\}$ is more similar to that of the original system.

Thus, we consider the
leader selection and demotion problems.
The leader selection problem is formulated as a minimization problem for the $H^2$ norm of the difference between the transfer functions of the original and new agent systems, under the assumption that the leader agents to be demoted are fixed. 
Here, the leaders of the new agent system are selected from all agents  except for demoted leader agents from the original system, and the number of new leaders is equal to the difference between the numbers of original and demoted leaders.
Because the problem is a finite combinatorial optimization problem, a brute force approach for solving the problem quickly becomes intractable as the number of leaders increases.
The leader demotion problem is that of finding optimal leader agents to be demoted, and is formulated
using the global optimal solution to the leader selection problem.

The contributions of this paper are as follows.\\
1) We prove that a global optimal solution to the leader selection problem is the set of original leader agents except for those that are demoted to followers.
To this end, we relax the original problem into a differentiable problem.
Then, by calculating the gradient and Hessian of the objective function of the relaxed problem, we prove that the function is convex.
It is shown that zero points of the gradient are global optimal solutions to the leader selection problem, which is a finite combinatorial optimization problem. 
Furthermore, we prove that any set of leader agents to be demoted subject to a fixed number of elements is a solution to the leader demotion problem.

\noindent
2) By combining the solutions to the leader selection and demotion problems, we prove that
 if we choose new leader agents from the original leader agents except for those specified by the set of leader agents to be demoted,
then the relative $H^2$ error between the  transfer functions of the original and new agent systems is completely determined by the numbers of original leader agents and  leader agents that are demoted to follower agents.
That is, we reveal that the relative $H^2$ error does not depend on the number of agents on the graph.

The remainder of this paper is organized as follows.
In Section \ref{sec2}, we formulate the leader selection and demotion problems.
In Section \ref{sec3}, we rigorously solve these problems.
In Section \ref{sec4}, we verify the solutions using a simple example. 
Finally, our conclusions are presented in Section \ref{sec5}.

{\it Notation:} The sets of real and complex numbers are denoted by ${\bf R}$ and ${\bf C}$, respectively.
The identity matrix of size $n$ is denoted by $I_n$.
The symbol ${\bf 0}_n \in {\bf R}^n$ is a vector with all zero entries.
For any finite set $S$, $|S|$ denotes the number of elements of $S$.
Given a vector $v\in {\bf C}^n$, $||v||$ denotes the Euclidean norm.
The Hilbert space $L^2({\bf R}^n)$ is defined by
\begin{align*}
L^2({\bf R}^n) := \left\{ f: [0,\infty) \rightarrow {\bf R}^n\, \big|\, \int_0^{\infty} ||f(t)||^2 dt <\infty \right\}.
\end{align*}
Given a measurable function $f: [0,\infty)\rightarrow {\bf R}^n$, $||f||_{L^2}$ and $||f||_{L^{\infty}}$ denote the $L^2$ and $L^{\infty}$ norms of $f$, respectively, i.e.,
\begin{align*}
& ||f||_{L^2}:= \sqrt{ \int_0^{\infty} ||f(t)||^2 dt }, \\
& ||f||_{L^{\infty}} := \sup_{t\geq 0} ||f(t)||. 
\end{align*}
Given a matrix $A\in {\bf C}^{m\times n}$, $||A||$ and $||A||_F$ denote the induced norm and the Frobenius norm, respectively, i.e.,
\begin{align*}
& ||A|| := \max_{v\in {\bf C}^n\backslash \{0\}} \frac{ ||Av|| }{||v||}, \\
& ||A||_F:= \sqrt{ {\rm tr} (A^{\dagger} A)},
\end{align*}
where the superscript $\dagger$ denotes the Hermitian conjugation, and
${\rm tr} (A^{\dagger}A)$ is the trace of $A^{\dagger}A$, i.e., the sum of the diagonal elements of $A^{\dagger}A$.
For a matrix function $G(s)\in {\bf C}^{m\times n}$, $||G||_{H^2}$ denotes the $H^2$ norm of $G$; i.e.,
\begin{align*}
||G||_{H^2}  := \sqrt{ \frac{1}{2\pi} \int_{-\infty}^{\infty} ||G(i\omega)||_F^2 d \omega},
\end{align*}
where $i$ is the imaginary unit.

\section{Problem formulation} \label{sec2}

\subsection{Leader-follower multi-agent systems} 

This subsection defines leader-follower multi-agent systems, to be studied in this paper.

Let $\mathcal{G}=(V,E,A)$ be a connected weighted undirected graph, where $V=\{1,2,\ldots, n\}$ is the vertex set, $E$ is the edge set, and $A$ is the adjacency matrix
consisting of nonnegative elements $a_{ij}$ called the weights.
That is, for each edge $(i,j)\in E$, the $i$-th row and $j$-th column entry of $A$ is equal to the weight $a_{ij}$, and all other entries of $A$ are equal to zero.
Because $\mathcal{G}$ is an undirected graph, $a_{ij}=a_{ji}$; i.e., the matrix $A$ is symmetric.
The degree matrix of $\mathcal{G}$ is a diagonal matrix denoted by $D={\rm diag}(d_1,d_2,\ldots, d_n)$, with
\begin{align*}
d_i = \sum_{j=1}^n a_{ij}.
\end{align*}
The Laplacian matrix of $\mathcal{G}$ is defined as $L=D-A$.
Let the total number of edges be $k$.
Then, we number the edges of $G$ by a unique $e\in \{1,2,\ldots,k\}$, and assign an arbitrary direction to each edge.
The incidence matrix $R=(R_{ij})\in {\bf R}^{n\times k}$ of the graph $G$ is defined by
\begin{align*}
R_{ij} := \begin{cases}
1,\quad\,\,\,\,\, {\rm if\, vertex\,}\, i\, {\rm is\, the\, source\, node\, of\, edge}\, j \\
-1,\quad {\rm if\, vertex\,}\, i\, {\rm is\, the\, sink\, node\, of\, edge}\, j \\
0,\quad\,\,\,\,\, {\rm otherwise}.
\end{cases}
\end{align*}
Furthermore, let
\begin{align*}
W:= {\rm diag} (w_1,w_2,\ldots,w_k)
\end{align*}
be the diagonal matrix of edge weights.
Then, the relationship between $L$, $R$, and $W$ is described by
\begin{align}
L = RWR^T, \label{L2}
\end{align}
which can be found in \cite{mesbahi2010graph}.

Let $V_L:= \{v_1,v_2,\ldots,v_m\}\subset V$ be a set of leader agents and $V_F:=V\backslash V_L$ a set of follower agents.
We consider the following leader-follower multi-agent system:
\begin{align}
\dot{x}_i = \begin{cases}
\sum_{j=1}^n a_{ij} (x_j-x_i)\quad\quad\quad\, {\rm if}\quad i\in V_F,\\
\sum_{j=1}^n a_{ij} (x_j-x_i)+u_l\quad {\rm if}\quad i=v_l\in V_L,
\end{cases} \label{multi_agent1}
\end{align}
where $x_i\in {\bf R}$ and $u_l\in {\bf R}$ denote the state of agent $i\in V$ and external input applied to agent $v_l\in V_L$, respectively. 
In \cite{monshizadeh2014projection}, the same multi-agent system has been considered.
The system \eqref{multi_agent1} can be rewritten as
\begin{align}
\dot{x} = -Lx + Mu, \label{sys}
\end{align}
where $x:=(x_1,x_2,\ldots, x_n) \in {\bf R}^n$, $u:=(u_1,u_2,\ldots, u_m) \in {\bf R}^m$, and 
\begin{align*}
M_{il} = \begin{cases}
1\quad {\rm if}\quad i=v_l\in V_L, \\
0\quad {\rm otherwise}.
\end{cases}
\end{align*}
Moreover, we choose the output variable $y\in {\bf R}^k$ as
\begin{align}
y = W^{1/2}R^T x. \label{output}
\end{align}
This is because the output $y$ defined by \eqref{output} reflects the disagreements among agents.
In fact, it follows from \eqref{L2} that $||y||^2 = x^T L x = \frac{1}{2} \sum_{i,j} a_{ij} (x_i-x_j)^2$, which gives a measure of group disagreement \cite{olfati2004consensus}.
Then, the transfer function of the input-output system \eqref{sys}-\eqref{output} is defined as
\begin{align}
G(s) := W^{1/2}R^T (sI_n+L)^{-1} M \label{transfer1}
\end{align}
for $s\in {\bf C}$.

\subsection{Leader selection and demotion problems}

This subsection formulates the leader selection and demotion problems to be solved in this paper.

Let $V_{\rm demotion}\subset V_L$ be the set of leader agents that are demoted to follower agents, 
$\tilde{V}_L\subset V\backslash V_{\rm demotion}$ the new set of leader agents, and $\tilde{V}_F:= V\backslash \tilde{V}_L$ the new set of follower agents.
Here, we assume that
\begin{align}
|\tilde{V}_L|= m -|V_{\rm demotion}|. \label{assumption}
\end{align}
That is, the number of new leader agents is equal to the difference between the number of original leader agents and the number of follower agents demoted from leader agents.  
If $|V_{\rm demotion}| = r$, then we can express $V_{\rm demotion}$ as  $\{ v_{i_1}, v_{i_2}, \ldots, v_{i_r}\}$, where $i_1, i_2,\ldots, i_r\in \{1,2,\ldots, m\}$.
Thus, from the assumption \eqref{assumption}, $\tilde{V}_L\subset V\backslash V_{\rm demotion}$ can be written as 
\begin{align*}
\tilde{V}_L =\{ \tilde{v}_1,\tilde{v}_2,\ldots, \tilde{v}_m\}\backslash \{ \tilde{v}_{i_1},\tilde{v}_{i_2},\ldots, \tilde{v}_{i_r}\},
\end{align*}
where $\{ v_{i_1}, v_{i_2}, \ldots, v_{i_r}\}=\{ \tilde{v}_{i_1},\tilde{v}_{i_2},\ldots, \tilde{v}_{i_r}\}$.
Note that $V_L\neq \{ \tilde{v}_1,\tilde{v}_2,\ldots, \tilde{v}_m\}$ in general.

In addition to the original leader-follower multi-agent system \eqref{multi_agent1}, we also consider the following system
\begin{align}
\dot{\tilde{x}}_i = \begin{cases}
\sum_{j=1}^n a_{ij} (\tilde{x}_j-\tilde{x}_i)\quad\quad\quad\, {\rm if}\quad i\in \tilde{V}_F,\\
\sum_{j=1}^n a_{ij} (\tilde{x}_j-\tilde{x}_i)+u_l\quad {\rm if}\quad i=\tilde{v}_l\in \tilde{V}_L.
\end{cases} \label{multi_agent2}
\end{align}
The system \eqref{multi_agent2} can be rewritten as
\begin{align}
\dot{\tilde{x}} = -L\tilde{x} + \tilde{M}u, \label{sys2}
\end{align}
where
\begin{align}
\tilde{M}_{il} = \begin{cases}
1\quad {\rm if}\quad i=\tilde{v}_l\in \tilde{V}_L, \\
0\quad {\rm otherwise}.
\end{cases} \label{tilde_M}
\end{align}
As is the case with the original system, we choose the output $\tilde{y}\in {\bf R}^k$ as 
\begin{align}
\tilde{y}=W^{1/2}R^T\tilde{x}. \label{output2}
\end{align}
Then, the transfer function of the input-output system \eqref{sys2}-\eqref{output2} is defined as 
\begin{align}
\tilde{G}(s) := W^{1/2}R^T (sI_n+L)^{-1} \tilde{M}. \label{transfer2}
\end{align}

Note that although there is a non-zero entry for each column of the matrix $M$, all the entries of certain columns of the matrix $\tilde{M}$, specified by the subscripts of the entries in $V_{\rm demotion}$ are zeros. 
For example, let $n=6$ and $m=3$. Assume that $V_L = \{ v_1,v_2, v_3\}$, with $v_1=1$, $v_2=4$, and $v_3=5$.
Then, the matrix $M$ of \eqref{sys} is given by
\begin{align*}
M = \begin{pmatrix}
1 & 0 & 0 \\
0 & 0 & 0 \\
0 & 0 & 0 \\
0 & 1 & 0 \\
0 & 0 & 1 \\
0 & 0 & 0
\end{pmatrix}.
\end{align*}
Thus, there is a non-zero entry for each column of $M$.
Next, let $V_{\rm demotion} = \{v_2\}=\{\tilde{v}_2\}$ and $\tilde{V}_L = \{\tilde{v}_1,\tilde{v}_2,\tilde{v}_3\}\backslash \{\tilde{v}_2\} = \{\tilde{v}_1,\tilde{v}_3\}$, with $\tilde{v}_1=2$ and $\tilde{v}_3=6$.
Then, the matrix $\tilde{M}$ of \eqref{sys2} is given by
\begin{align}
\tilde{M} = \begin{pmatrix}
0 & 0 & 0 \\
1 & 0 & 0 \\
0 & 0 & 0 \\
0 & 0 & 0 \\
0 & 0 & 0 \\
0 & 0 & 1
\end{pmatrix}. \label{ex}
\end{align}
Thus, all entries of the second column, specified by the subscript of the entry in $V_{\rm demotion}=\{\tilde{v}_2 \}$, are zeros.

In this paper, we want to determine the set $\tilde{V}_L\subset V\backslash V_{\rm demotion}$ minimizing 
\begin{align*}
f(\tilde{V}_L) := ||G-\tilde{G}||_{H^2}^2.
\end{align*}
This is because the following lemma holds.
\begin{lemma} \label{Lem1}
If $u\in L^2({\bf R}^m)$, then 
\begin{align}
||y-\tilde{y}||_{L^{\infty}} \leq ||G-\tilde{G}||_{H^2}\cdot ||u||_{L^2}.
\end{align}
\end{lemma}
This means that if $f(\tilde{V}_L)$ and $||u||_{L^2}$ are sufficiently small, then $||y(t)-\tilde{y}(t)||$  is also small for any $t\geq 0$.
Although in \cite{gugercin2008h_2} and \cite{sato2017riemannian} similar results to Lemma \ref{Lem1} were proved for asymptotically stable systems,
the systems \eqref{sys} and \eqref{sys2} are not asymptotically stable because the matrix $-L$ has at least one zero eigenvalue.
We provide a proof of Lemma \ref{Lem1} in Appendix \ref{proof_Lem1}.

To this end, we solve the following leader selection problem.

\begin{framed}
Problem 1 (Leader selection problem): 
\begin{align*}
&{\rm Given}\,\, V_{\rm demotion}\subset V_L,\\
& {\rm find}\,\, \tilde{V}_L\subset V\backslash V_{\rm demotion} \,\, 
{\rm minimizing}\,\, f(\tilde{V}_L) \\
&{\rm subject\, to}\,\, \eqref{assumption}.
\end{align*}
\end{framed}

\noindent
By solving Problem 1, we can find new leader agents when the original leader agents that are demoted to follower agents have been determined, in the sense
that the transfer function \eqref{transfer2} best approximates \eqref{transfer1} in the sense of the $H^2$ norm.

However, if the number of agents $n$ is large, then it is difficult to solve Problem 1 by 
 enumerating all possible subsets of size $\binom{n -|V_{\rm demotion}|}{m -|V_{\rm demotion}|}$, evaluating $f$ for all of these subsets,
and picking the best subset.
In the next section, we provide a global optimal solution to Problem 1.

Because a solution to Problem 1 may depend on $V_{\rm demotion}$,
we want to determine the set $V_{\rm demotion}\subset V_L$ minimizing
\begin{align}
g(V_{\rm demotion}) := f( \tilde{V}^*_L ), \label{g_def}
\end{align}
where $\tilde{V}^*_L$ is a global optimal solution to Problem 1.
To this end,
we further consider the following leader demotion problem.

\begin{framed}
Problem 2 (Leader demotion problem): 
\begin{align*}
&{\rm Given}\,\, r\in \{1,2,\ldots,m\},\\
& {\rm find} \,\, V_{\rm demotion}\,\, {\rm minimizing}\,\, g(V_{\rm demotion})\\
&{\rm subject\, to}\,\, |V_{\rm demotion}|=r.
\end{align*}
\end{framed}

\noindent
In Section \ref{Sec3B}, we solve Problem 2.


\section{Solutions to Problems 1 and 2} \label{sec3}

\subsection{Solution to Problem 1}

This subsection proves the following theorem.

\begin{theorem} \label{thm}
A global optimal solution
$\tilde{V}_L^*$
 to Problem 1 is given by
\begin{align}
\tilde{V}_L^* = V_L\backslash V_{\rm demotion}. \label{result1}
\end{align}
\end{theorem}

This theorem means that
new leader agents minimizing $||G-\tilde{G}||_{H^2}$ are composed of the original leader agents except for those that are demoted to follower agents.
Note that this theorem holds for any connected weighted undirected graph $\mathcal{G}$.

To prove this theorem, we note that there exists a function $h$ satisfying
\begin{align}
h(\tilde{M}) = f(\tilde{V}_L), \label{touka}
\end{align}
because there is a one-to-one relation between the matrix $\tilde{M}$ and the set $\tilde{V}_L$.
Furthermore, note that
\begin{align*}
\tilde{M} \in Z_J^{n\times m} \subset {\bf R}^{n\times m}_J \subset {\bf R}^{n\times m}.
\end{align*}
Here, we assume the following:
\begin{itemize}
\item $J:=\{i_1,i_2,\ldots, i_r\}$, which constitutes the set of subscripts of the elements of $V_{\rm demotion}=\{ v_{i_1}, v_{i_2},\ldots, v_{i_r}\}$.
\item $Z_J^{n\times m}$ is the set of $n$-row and $m$-column matrices composed of $0$ and $1$. The $i_1,i_2,\ldots, i_r$-th column vectors of matrices in $Z_J^{n\times m}$ are zero vectors.
The other column vectors are non-zero vectors, but consist of zeros except for one entry. Furthermore, the other column vectors are not equal to each other.
For example, the matrix $\tilde{M}$ in \eqref{ex} is contained in $Z_{\{2\}}^{6\times 3}$.
\item ${\bf R}^{n\times m}_J$ is the set of free real matrices except for the $i_1,i_2,\ldots, i_r$-th column vectors. The $i_1,i_2,\ldots, i_r$-th column vectors of matrices in ${\bf R}^{n\times m}_J$ are zero vectors. That is, the set ${\bf R}^{n\times m}_J$ is a subspace of ${\bf R}^{n\times m}$.
\end{itemize}

From \eqref{touka},
Problem 1 is equivalent to the following problem.
\begin{framed}
Problem 3:
\begin{align*}
&{\rm Given}\,\, J\subset \{1,2,\ldots, m\},\\
& {\rm find}\,\, \tilde{M}\in Z_J^{n\times m} \,\, 
{\rm minimizing}\,\, h(\tilde{M}).
\end{align*}
\end{framed}
Problem 3 is relaxed to the following problem, because $Z_J^{n\times m}\subset {\bf R}^{n\times m}_J$.
\begin{framed}
Problem 4:
\begin{align*}
&{\rm Given}\,\, J\subset \{1,2,\ldots, m\},\\
& {\rm find}\,\, \tilde{M}\in {\bf R}_J^{n\times m} \,\, 
{\rm minimizing}\,\, h(\tilde{M}).
\end{align*}
\end{framed}
Let $\tilde{M}^*\in Z_J^{n\times m}$ be a global optimal solution to Problem 3. Because Problem 4 is a relaxation of Problem 3, there exists a solution $\tilde{M}\in {\bf R}_J^{n\times m}$ to Problem 4 such that
\begin{align*}
h(\tilde{M}) \leq h(\tilde{M}^*).
\end{align*}

In what follows, we show that
the global optimal solution $\tilde{M}\in {\bf R}_J^{n\times m}$ to Problem 4 coincides with $\tilde{M}^*\in Z_J^{n\times m}$, and
\begin{align}
\tilde{M}^*= M_J, \label{result2}
\end{align}
where $M_J$ denotes the same matrix as $M$ except for the $i_1,i_2,\ldots, i_r$-th column vectors, and each of these column vectors is replaced by a zero vector.
Because $\tilde{M}$ is defined by \eqref{tilde_M}, \eqref{result2} is equivalent to \eqref{result1}.

To proceed, we use the fact that
\begin{align}
h(\tilde{M}) = {\rm tr} ( (M-\tilde{M})^T W_o (M-\tilde{M}) ),  \label{key1}
\end{align}
where $W_o$ is given by \eqref{key2} in Appendix \ref{proof_Lem1}.
Note that $W_o$ is the observability Gramian of the system \eqref{sys}-\eqref{output} \cite{zhou1996robust}.
Eq.\,\eqref{key1} follows from a similar discussion as the derivation of \eqref{key1proof} and Parseval's theorem.

Now, we show that the objective function $h$ is convex on ${\bf R}^{n\times m}_J$ by calculating the gradient and Hessian.
To this end, let $\bar{h}$ denote the extension of the objective function $h$ to the ambient Euclidean space ${\bf R}^{n\times m}$.
The directional derivative of $\bar{h}$ at $\tilde{M}$ in the direction $\tilde{M}'$ can be calculated as
\begin{align}
{\rm D}\bar{h}(\tilde{M})[\tilde{M}']  = {\rm tr}(\tilde{M}'^T (-2W_o(M-\tilde{M}))). \label{1000}
\end{align}
Here, the directional derivative is defined in Appendix \ref{frechet}.
Because the gradient $\nabla \bar{h}(\tilde{M})$ satisfies
${\rm D}\bar{h}(\tilde{M})[\tilde{M}'] =  {\rm tr} (\tilde{M}'^T (\nabla \bar{h}(\tilde{M})))$,
(\ref{1000}) implies that
\begin{align*}
 \nabla \bar{h}(\tilde{M})  = -2W_o(M-\tilde{M}). 
\end{align*}

\noindent
Because $(\nabla \bar{h}(\tilde{M}))_{J}$ is the projection of  
$\nabla \bar{h}(\tilde{M})$ onto the subspace ${\bf R}^{n\times m}_J$,
the gradient ${\rm grad}\, h(\tilde{M})$ for $\tilde{M}\in {\bf R}^{n\times m}_J$ is given by
\begin{align}
{\rm grad}\, h(\tilde{M}) =& (\nabla \bar{h}(\tilde{M}))_{J} \nonumber \\
=& -2W_o (M-\tilde{M})_{J}\nonumber \\
=& -2W_o\left( M_{J}-\tilde{M} \right). \label{gradient}
\end{align}
The Hessian ${\rm Hess}\, h(\tilde{M})$ at any $\tilde{M}\in {\bf R}^{n\times m}_J$ is given by
\begin{align}
{\rm Hess}\, h(\tilde{M})[\tilde{M}'] =& \left( {\rm D} {\rm grad} h(\tilde{M}) [\tilde{M}'] \right)_{J} \nonumber\\
=& 2W_o \tilde{M}', \label{Hess}
\end{align}
where $\tilde{M}'\in T_{\tilde{M}} {\bf R}^{n\times m}_J$, and $T_{\tilde{M}} {\bf R}^{n\times m}_J$ is the tangent space of ${\bf R}^{n\times m}_J$ at the point $\tilde{M}$.
Here, note that $T_{\tilde{M}} {\bf R}^{n\times m}_J$ can be identified with ${\bf R}^{n\times m}_J$, because ${\bf R}^{n\times m}_J$ is a vector space.
For a detailed explanation of the concept of the Hessian, see \cite{absil2009optimization}.
Thus,
\begin{align}
\langle \tilde{M}', {\rm Hess}\, h(\tilde{M}) [\tilde{M}'] \rangle :=& {\rm tr} (\tilde{M}'^T{\rm Hess}\, h(\tilde{M}) [\tilde{M}'] ) \nonumber \\
=& 2{\rm tr} (\tilde{M}'^T W_o \tilde{M}'). \label{Hess2}
\end{align}
Because the observability Gramian $W_o$ is symmetric positive semidefinite, (\ref{Hess2}) implies that
$\langle \tilde{M}', {\rm Hess}\, h(\tilde{M}) [\tilde{M}'] \rangle \geq 0$ for any $0\neq \tilde{M}'\in T_{\tilde{M}} {\bf R}^{n\times m}_J$ and any $\tilde{M}\in {\bf R}^{n\times m}_J$.
Hence, the objective function $h$ is convex on ${\bf R}^{n\times m}_J$ \cite{boyd2004convex}.

If 
\begin{align}
\tilde{M} = M_J, \label{main}
\end{align}
then (\ref{gradient}) yields that ${\rm grad}\,h(\tilde{M})=0$; i.e.,
(\ref{main}) is at least a local optimal solution.
In fact, (\ref{main}) is a global optimal solution to Problem 4, because the function $h$ is  convex on ${\bf R}^{n\times m}_J$.
Because $M_J \in Z_J^{n\times m}$, \eqref{main} is also a global optimal solution to Problem 3.
As mentioned previously, \eqref{main} is equivalent to \eqref{result1}. 
This completes the proof.

\subsection{Solution to Problem 2}  \label{Sec3B}

This subsection proves the following theorem.

\begin{theorem}
Any $V_{\rm demotion}\subset V_L$ subject to $|V_{\rm demotion}|=r$ is a solution to Problem 2, and
\begin{align}
g(V_{\rm demotion}) = \frac{r}{2} \left(1-\frac{1}{n}\right). \label{result3}
\end{align}
\end{theorem}

Note that this theorem holds for any connected weighted undirected graph $\mathcal{G}$.

As a corollary of Theorems 1 and 2, we obtain the following.

\begin{corollary}
Suppose that $V_{\rm demotion}\subset V_L$ subject to $|V_{\rm demotion}|=r$, and $\tilde{V}^*_L$ is given by \eqref{result1}.
Then, 
\begin{align}
\frac{||G-\tilde{G}||_{H^2}}{||G||_{H^2}} = \sqrt{\frac{r}{m}}. \label{result4}
\end{align}
\end{corollary}

This corollary means that if we choose new leader agents from the original ones except for those specified by $V_{\rm demotion}$, then
the relative $H^2$ error between the transfer functions $G$ and $\tilde{G}$ is completely determined by $m$, the number of original leader agents, and $r$, the number of leader agents that are demoted to follower agents. That is, the relative $H^2$ error does not depend on the number of agents $n$ on the graph $\mathcal{G}$. 
Note that this corollary also holds for any connected weighted undirected graph $\mathcal{G}$.

First, we prove Theorem 2.
From \eqref{g_def}, \eqref{touka}, and \eqref{key1}, we obtain that
\begin{align}
g(V_{\rm demotion}) &= f(\tilde{V}_L^*) \nonumber\\
 &= {\rm tr} ( (M-\tilde{M}^*)(M-\tilde{M}^*)^T W_o), \label{key3}
\end{align}
where $\tilde{M}^*$ is given by \eqref{result2}.
Here, the second equality follows from the property that the trace is invariant under cyclic permutations.
The matrix $(M-\tilde{M}^*)(M-\tilde{M}^*)^T$ in \eqref{key3} is a diagonal matrix, where the diagonal elements are either zero or one.
Furthermore, the $i$-th diagonal element equals one if $i\in V_{\rm demotion}$ and zero otherwise.
Hence, it follows from \eqref{key2} and \eqref{key3} that \eqref{result3} holds.

Next, we prove Corollary 1.
Through the same discussion as for the derivations of \eqref{key1} and \eqref{key3}, we obtain
\begin{align}
||G||_{H^2}^2 = {\rm tr} (M^TW_oM) = {\rm tr} (MM^TW_o). \label{key4}
\end{align}
The matrix $MM^T$ in \eqref{key4} is also a diagonal matrix, where the diagonal elements are either zero or one.
Furthermore, the $i$-th diagonal element equals one if $i\in V_{L}$ and zero otherwise.
Thus, it follows from \eqref{key2} and \eqref{key4} that 
\begin{align*}
||G||_{H^2}^2 = \frac{m}{2} \left(1-\frac{1}{n} \right).
\end{align*}
From this and \eqref{result3}, we obtain \eqref{result4}.

\section{Numerical example} \label{sec4}

This section verifies Theorems 1 and 2 using a simple example.

Let $n=5$, $m=3$, $V=\{1,2,3,4,5\}$, $V_L=\{1,2,3\}$, and
\begin{align*}
A = \begin{pmatrix}
0 & 1 & 1 & 0  & 0 \\
1 & 0 & 0 & 1 & 0 \\
1 & 0 & 0 & 1 & 1 \\
0 & 1 & 1 & 0 & 1 \\
0 & 0 & 1 & 1 & 0
\end{pmatrix}.
\end{align*}

Table \ref{table1} shows the relation between $\tilde{V}_L\subset V\backslash V_{\rm demotion} = \{2,3,4,5\}$, i.e., $V_{\rm demotion}=\{1\}$,  and $f(\tilde{V}_L)$. 
Table \ref{table2} shows the relation between $\tilde{V}_L\subset V\backslash V_{\rm demotion} = \{1,3,4,5\}$, i.e., $V_{\rm demotion}=\{2\}$,  and $f(\tilde{V}_L)$.
Table \ref{table3} shows the relation between $\tilde{V}_L\subset V\backslash V_{\rm demotion} = \{1,2,4,5\}$, i.e., $V_{\rm demotion}=\{3\}$,  and $f(\tilde{V}_L)$.

As shown in Tables \ref{table1}, \ref{table2}, and \ref{table3}, $f(\tilde{V}_L)$ is minimized if \eqref{result1} in Theorem 1 holds.
Furthermore, we can observe in Tables \ref{table1}, \ref{table2}, and \ref{table3} that for any $V_{\rm demotion} \subset V_L$ subject to $|V_{\rm demotion}| =1$,
$g(V_{\rm demotion})$ is minimized and $g(V_{\rm demotion})=0.4000= \frac{1}{2}(1-\frac{1}{5})$; i.e., \eqref{result3} in Theorem 2 holds.

\begin{table}[t]
\caption{The relation between $\tilde{V}_L\subset \{2,3,4,5\}$ and $f(\tilde{V}_L)$.} \label{table1}
  \begin{center}
    \begin{tabular}{|c|c|c|c|c|c|c|} \hline
      $\tilde{V}_L$    &   $\{2,3\}$   &  $\{2,4\}$    & $\{2,5\}$ &  $\{3,4\}$     & $\{3,5\}$      & $\{4,5\}$  \\ \hline 
      $f(\tilde{V}_L)$ &   0.4000      & 1.4000     &     1.4000  & 2.4000 & 2.4000 & 2.4000 \\\hline 
    \end{tabular}
  \end{center}
\end{table}

\begin{table}[t]
\caption{The relation between $\tilde{V}_L\subset \{1,3,4,5\}$ and $f(\tilde{V}_L)$.} \label{table2}
  \begin{center}
    \begin{tabular}{|c|c|c|c|c|c|c|} \hline
      $\tilde{V}_L$    &   $\{1,3\}$   &  $\{1,4\}$    & $\{1,5\}$ &  $\{3,4\}$     & $\{3,5\}$      & $\{4,5\}$  \\ \hline 
      $f(\tilde{V}_L)$ &   0.4000      & 1.4000     &     1.4000  & 2.4000 & 2.4000 & 2.4000 \\\hline 
    \end{tabular}
  \end{center}
\end{table}

\begin{table}[t]
\caption{The relation between $\tilde{V}_L\subset \{1,2,4,5\}$ and $f(\tilde{V}_L)$.} \label{table3}
  \begin{center}
    \begin{tabular}{|c|c|c|c|c|c|c|} \hline
      $\tilde{V}_L$    &   $\{1,2\}$   &  $\{1,4\}$    & $\{1,5\}$ &  $\{2,4\}$     & $\{2,5\}$      & $\{4,5\}$  \\ \hline 
      $f(\tilde{V}_L)$ &   0.4000      & 1.4000     &     1.4000  & 2.4000 & 2.4000 & 2.4000 \\\hline 
    \end{tabular}
  \end{center}
\end{table}

\section{Conclusion} \label{sec5}

We have considered the leader selection and demotion problems.
We have proved that a global optimal solution to the leader selection problem is the set of original leader agents except for those that are demoted to followers.
Furthermore, we have proved that any set of leader agents to be demoted subject to a fixed number of elements is a solution to the leader demotion problem.
By combining the solutions to the leader selection and demotion problems, we have also proved that
 if we choose new leader agents from the original leader agents except for those specified by the set of leader agents to be demoted, then
the relative $H^2$ error between the  transfer functions of the original and new agent systems is completely determined by the numbers of original leader agents and  leader agents that are demoted to follower agents.
That is, we have revealed that the relative $H^2$ error does not depend on the number of agents on a connected weighted undirected graph. 
Finally, we have verified the solutions using a simple example.

\appendix
\subsection{Proof of Lemma \ref{Lem1}} \label{proof_Lem1}

We first prove that the impulse response of the system \eqref{sys}-\eqref{output}
\begin{align*}
g(t) = \begin{cases}
0\quad \quad\quad\quad\quad\quad\quad\quad\quad\quad\quad (t<0),\\
W^{1/2}R^T\exp (-Lt) M\quad\quad\,\, (t\geq 0)
\end{cases}
\end{align*}
 is contained in the $L^2$ space given by the Hilbert space of matrix-valued functions on ${\bf R}$ with
the inner product $\langle f_1, f_2 \rangle := \int_{-\infty}^{\infty} {\rm tr} (f_1(t)^Tf_2(t)) dt$.
When $t\geq 0$, we obtain
\begin{align*}
||g(t)||_F^2 = {\rm tr} ( M^T \exp(-Lt)L \exp(-Lt) M).
\end{align*}
Thus,
\begin{align}
\int_{-\infty}^{\infty} ||g(t)||_F^2 dt &= \int_0^{\infty} ||g(t)||_F^2 dt \nonumber\\
&= {\rm tr}( M^T W_o M ), \label{key1proof}
\end{align}
where
\begin{align}
W_o &= \int_0^{\infty} \exp (-Lt) L \exp (-Lt) dt \nonumber\\
&= \frac{1}{2} I_n - \frac{1}{2n} {\bf 1} {\bf 1}^T. \label{key2}
\end{align}
These expressions can be also found in the proof of Theorem 6 in \cite{monshizadeh2014projection}.
Hence, the impulse response $g$ is contained in the $L^2$ space.

Because $g$ is contained in the $L^2$ space, the Fourier transformation of $g$ can be defined \cite{zhou1996robust}; i.e., $G(i\omega)$ can be defined.
From the same discussion, $\tilde{G}(i\omega)$ can also be defined.
Thus, if $u\in L^2({\bf R}^m)$, then $Y(i\omega) =G(i\omega) U(i\omega)$ and $\tilde{Y}(i\omega)=\tilde{G}(i\omega)U(i\omega)$ can be defined,
where $U$, $Y$, and $\tilde{Y}$ are the Fourier transformations of $u$, $y$, and $\tilde{y}$, respectively.
Hence, if $u\in L^2({\bf R}^m)$, then we have that
\begin{align*}
& ||y-\tilde{y}||_{L^{\infty}}\\
=& \sup_{t\geq 0} ||y(t)-\tilde{y}(t)|| \\
=& \sup_{t\geq 0} ||\frac{1}{2\pi} \int_{-\infty}^{\infty} (Y(i\omega)-\tilde{Y}(i\omega)) e^{i\omega t} d\omega || \\
\leq & \frac{1}{2\pi} \int_{-\infty}^{\infty} ||Y(i\omega)|| d\omega \\
\leq & \frac{1}{2\pi} \int_{-\infty}^{\infty} ||G(i\omega)-\tilde{G}(i\omega)|| \cdot ||U(i\omega)|| d\omega \\
\leq & \sqrt{\frac{1}{2\pi} \int_{-\infty}^{\infty} ||G(i\omega)-\tilde{G}(i\omega)||^2d\omega}\sqrt{\frac{1}{2\pi} \int_{-\infty}^{\infty} ||U(i\omega)||^2d\omega}\\
\leq &||G-\tilde{G}||_{H^2}\cdot ||u||_{L^2},
\end{align*}
where the second equality follows from the inverse Fourier transformations of $Y$ and $\tilde{Y}$, the fifth inequality is from the Cauchy-Schwarz inequality, and the
final inequality follows from $||G(i\omega)-\tilde{G}(i\omega)||\leq ||G(i\omega)-\tilde{G}(i\omega)||_F$ and Parseval's theorem.
This completes the proof.

\subsection{Directional derivative of smooth functions} \label{frechet}

Let $h$ and $\langle \cdot, \cdot \rangle$ be a smooth real-valued function on a finite-dimensional Euclidean space $E$ and the Euclidean inner product on $E$, respectively.
The Fr\'echet derivative ${\rm D}h(p):E\rightarrow {\bf R}$ of $h$ at $p\in E$ is defined as a linear operator such that
\begin{align*}
\lim_{\xi\rightarrow 0} \frac{ || h(p+\xi)-h(p)-{\rm D}h(p)[\xi]||}{||\xi||} = 0,
\end{align*}
where $||\cdot||$ is the Euclidean norm \cite{absil2009optimization}.
Then, ${\rm D}h(p)[\xi]$ is the directional derivative of $h$ at $p$ along $\xi\in E$ and the Euclidean gradient $\nabla h(p)$ at $p\in E$ satisfies
\begin{align*}
{\rm D}h(p)[\xi] = \langle \nabla h(p), \xi \rangle.
\end{align*}


%



\ifCLASSOPTIONcaptionsoff
  \newpage
\fi

\end{document}